\documentclass[a4paper,11pt]{article}
\usepackage{amsfonts,amssymb,amsmath}

\usepackage[english]{babel}
 \makeatletter
 \@addtoreset{equation}{section}
 \makeatother

 \makeatletter
 \@addtoreset{enunciato}{section}
 \makeatother

 \newcounter{enunciato}[section]

 \newtheorem{ittheorem}{Theorem}
 \newtheorem{itlemma}{Lemma}
 \newtheorem{itproposition}{Proposition}
 \newtheorem{itcorollary}{Corollary}
 \newtheorem{itdefinition}{Definition}
 \newtheorem{itremark}{Remark}
 \newtheorem{itclaim}{Claim}
 \newtheorem{itfact}{Fact}
 \newtheorem{itconjecture}{Conjecture}

 \newenvironment{theorem}{\addtocounter{enunciato}{1}
 \begin{ittheorem}}{\end{ittheorem}}

 \newenvironment{lemma}{\addtocounter{enunciato}{1}
 \begin{itlemma}}{\end{itlemma}}

 \newenvironment{proposition}{\addtocounter{enunciato}{1}
 \begin{itproposition}}{\end{itproposition}}

 \newenvironment{corollary}{\addtocounter{enunciato}{1}
 \begin{itcorollary}}{\end{itcorollary}}

 \newenvironment{definition}{\addtocounter{enunciato}{1}
 \begin{itdefinition}}{\end{itdefinition}}

 \newenvironment{remark}{\addtocounter{enunciato}{1}
 \begin{itremark}}{\end{itremark}}

 \newenvironment{claim}{\addtocounter{enunciato}{1}
 \begin{itclaim}}{\end{itclaim}}

 \newenvironment{fact}{\addtocounter{enunciato}{1}
 \begin{itfact}}{\end{itfact}}

 \newenvironment{conjecture}{\addtocounter{enunciato}{1}
 \begin{itconjecture}}{\end{itconjecture}}

 \newcommand{\be}[1]{\begin{equation}\label{#1}}
 \newcommand{\ee}{\end{equation}}

 \newcommand{\bl}[1]{\begin{lemma}\label{#1}}
 \newcommand{\el}{\end{lemma}}

 \newcommand{\br}[1]{\begin{remark}\label{#1}}
 \newcommand{\er}{\end{remark}}

 \newcommand{\bt}[1]{\begin{theorem}\label{#1}}
 \newcommand{\et}{\end{theorem}}

 \newcommand{\bd}[1]{\begin{definition}\label{#1}}
 \newcommand{\ed}{\end{definition}}

 \newcommand{\bcl}[1]{\begin{claim}\label{#1}}
 \newcommand{\ecl}{\end{claim}}

 \newcommand{\bfact}[1]{\begin{fact}\label{#1}}
 \newcommand{\efact}{\end{fact}}

 \newcommand{\bp}[1]{\begin{proposition}\label{#1}}
 \newcommand{\ep}{\end{proposition}}

 \newcommand{\bc}[1]{\begin{corollary}\label{#1}}
 \newcommand{\ec}{\end{corollary}}

 \newcommand{\bcj}[1]{\begin{conjecture}\label{#1}}
 \newcommand{\ecj}{\end{conjecture}}

 \newcommand{\bpr}{\begin{proof}}
 \newcommand{\epr}{\end{proof}}

 \newcommand{\bprlem}[1]{\begin{proofof}{\it Lemma \ref{#1}}.\,\,}
 \newcommand{\eprlem}{\end{proofof}}

 \newcommand{\bprthm}[1]{\begin{proofof}{\it Theorem \ref{#1}}.\,\,}
 \newcommand{\eprthm}{\end{proofof}}

 \newcommand{\bi}{\begin{itemize}}
 \newcommand{\ei}{\end{itemize}}

 \newcommand{\ben}{\begin{enumerate}}
 \newcommand{\een}{\end{enumerate}}

 \newenvironment{proof}{\noindent {\em Proof}.\,\,}{\hspace*{\fill}$\halmos$\medskip}
 \newenvironment{proofof}{\noindent {\em Proof of\,\,}}{\hspace*{\fill}$\halmos$\medskip}
 \newcommand{\halmos}{\rule{1ex}{1.4ex}}

 \newcommand{\one}{{\mathchoice {1\mskip-4mu\mathrm l}
         {1\mskip-4mu\mathrm l}
         {1\mskip-4.5mu\mathrm l}
         {1\mskip-5mu\mathrm l}}}

\def \E {{\mathbb E}}

\def \N {{\mathbb N}}
\def \P {{\mathbb P}}

\def \R {{\mathbb R}}
\def \Z {{\mathbb Z}}

\def \ra {\rightarrow}
\def \ba {\begin{array}}
\def \ea {\end{array}}

\def \var {{\rm var}}
\def \Var {{\rm Var}}

\def \cE {{\mathcal E}}
\def \cF {{\mathcal F}}
\def \cG {{\mathcal G}}

\def \cN {{\mathcal N}}

\def\one{\rlap{\mbox{\small\rm 1}}\kern.15em 1}

\def\embf#1{\emph{\bf #1}}

\begin{document}
\title{A functional central limit theorem for regenerative chains}

\author{\renewcommand{\thefootnote}{\arabic{footnote}}
G.\ Maillard
\footnotemark[1]
\\
\renewcommand{\thefootnote}{\arabic{footnote}}
S.\ Sch\"opfer
\footnotemark[1]
}

\footnotetext[1]{
Institut de Math\'ematiques, Station 8,
\'Ecole Polytechnique F\'ed\'erale de Lausanne,
CH-1015 Lausanne, Switzerland,
{\sl gregory.maillard@epfl.ch},
{\sl samuel.schoepfer@epfl.ch}
}
\date{}
\maketitle

\begin{abstract} Using the regenerative scheme of \cite{comferfer02}, we establish a functional
central limit theorem (FCLT) for discrete time stochastic processes (chains) with summable memory
decay. Furthermore, under stronger assumptions on the memory decay, we identify the limiting variance
in terms of the process only.
As applications, we define classes of binary autoregressive processes and power-law Ising chains
for which the FCLT is fulfilled.

{\it Acknowledgment.} The research of GM is partially supported by the SNSF, grant $\#\, 200020-115964/1$.
We thank Roberto Fern\'andez for pointing out this problem and Thomas Mountford for valuable
conversations.
\end{abstract}



\section{Introduction and preliminaries}
\label{S1}


\subsection{Introduction}
\label{S1.1}

Chains are discrete-time stochastic processes with infinite memory that are natural extensions
of Markov chains when the associated process depends on its whole past. Such processes have been
extensively studied (see e.g.\ Fern\'andez and Maillard \cite{fermai05} and references therein),
but surprisingly very few is known about limit theorems. In this paper, we partially fill this gap
by establishing a functional central limit theorem.

Historically, the first central limit theorems for chains have been established, under strong ergodic
assumptions, by Cohn (1966) \cite{coh66}, Ibragimov and Linnik (1971) \cite{ibrlin71} (see also
Iosifescu and Grigorescu (1990) \cite{griios90}). More recently, the empirical entropies of chains
with exponential memory decay have been studied both in terms of their limit behavior (Gabrielli,
Galves and Guiol (2003) \cite{gabgalgui03}) and their large deviations (Chazottes and Gabrielli
(2005) \cite{chagab05}).

Limit theorems such as LCLT and LIL have been broadly studied for Markovian chains (see e.g.\ Meyn
and Tweedie \cite{meytwe93}). Regeneration methods, introduced by Chung (1967) \cite{chu67} and
refined by Chen (1999) \cite{che99}, have been used to divide the Markov chain into independent
random blocks in order to derive such limit theorems. This constitute a motivating challenge to
extend those techniques to the non-Markovian case.

In that paper, we prove a general FCLT (Theorem \ref{th1}) for chains satisfying the regenerative
scheme introduced in Comets, Fern\'andez and Ferrari \cite{comferfer02}. We give an explicit
expression for the associated limiting variance depending both on the original and regenerative
processes. As a corollary (Corollary \ref{cor1}), we give a more tractable condition on the memory
decay of the chain under which the FCLT is fulfilled. We also give a regime of the memory decay for
which the limiting variance can be expressed in terms on the original process only. As applications,
we derive FCLTs for autoregressive binary processes and Ising chains (Propositions \ref{propAutoregress}
and \ref{propIsing}).

The paper is organized as follows. In the rest of this section, we give some definitions
and preliminaries. In Section \ref{S2}, we state the main results. In Section \ref{S3}
we introduce binary autoregressive processes and power-law Ising chains for which we give
central limit theorems. Finally, Section \ref{S4} is devoted to the proofs.


\subsection{Notation and preliminary definitions}
\label{S1.2}

We consider a measurable space $(E,\cE)$ where $E$ is a finite alphabet and
$\cE$ is the discrete $\sigma$-algebra. We denote $(\Omega,\cF)$ the associated
product measurable space with $\Omega=E^{\Z}$. For each $\Lambda\subset\Z$ we
denote $\Omega_\Lambda = E^{\Lambda}$ and $\sigma_\Lambda$ for the restriction
of a configuration $\sigma\in\Omega$ to $\Omega_\Lambda$, namely the family
$(\sigma_i)_{i\in\Lambda}\in E^{\Lambda}$. Also, $\cF_\Lambda$ will denote the
sub-$\sigma$-algebra of $\mathcal{F}$ generated by cylinders based on $\Lambda$
($\cF_\Lambda$-measurable functions are insensitive to configuration values outside
$\Lambda$). When $\Lambda$ is an interval, $\Lambda=[k,n]$ with $k,n\in\Z$ such
that $k\le n$, we use the notation:
$\omega_{k}^{n}=\omega_{[k,n]}=\omega_{k},\ldots
,\omega_{n}$, $\Omega_k^n=\Omega_{[k,n]}$ and $\cF_k^n=\cF_{[k,n]}$. For
semi-intervals we denote also $\cF_{\le n}=\cF_{(-\infty,n]}$, etc. The
concatenation notation $\omega_\Lambda\,\sigma_\Delta$, where
$\Lambda\cap\Delta=\emptyset$, indicates the configuration on $\Lambda\cup\Delta$
coinciding with $\omega_i$ for $i\in\Lambda$ and with $\sigma_i$ for $i\in\Delta$.


\subsection{Chains}
\label{S1.3}

We start by briefly reviewing the well-known notions of chains in a shift-invariant setting.
In this particular case, chains are also called $g$-measures (see \cite{kea72}).

\bd{lis1}
A \embf{$g$-function} $g$ is a probability kernel
$g\colon\Omega_0\times\Omega_{-\infty}^{-1}\to [0,1]$, i.e.,
\be{g-funct}
\sum_{\omega_0\in\Omega_0}g\big(\omega_0\mid\omega_{-\infty}^{-1}\big)=1,
\qquad\omega_{-\infty}^{-1}\in\Omega_{-\infty}^{-1}.
\ee
The $g$-function $g$ is:
\begin{itemize}
\item[\rm{(i)}]
\embf{Continuous} if the function
$g\left(\omega_0 \mid \cdot\,\right)$ is continuous for each $\omega_0\in\Omega_0$, i.e.,
for all $\epsilon>0$, there exists $n\geq 0$ so that
\be{contlis}
\big|g\big(\omega_0\mid \omega_{-\infty}^{-1}\big)
-g\big(\sigma_0\mid \sigma_{-\infty}^{-1}\big)\big|<\epsilon
\ee
for all $\omega_{-\infty}^{0},\sigma_{-\infty}^{0}\in\Omega_{-\infty}^{0}$ with
$\omega_{-n}^{0}=\sigma_{-n}^{0}$;
\item[\rm{(ii)}]
\embf{Bounded away form zero} if $g\left(\omega_0 \mid \cdot\,\right)\geq c>0$ for each $\omega_0\in\Omega_0$;
\item[\rm{(iii)}]
\embf{Regular} if $f$ is continuous and bounded away from zero.
\end{itemize}
\ed
\bd{chain}
A probability measure $\P$ on $(\Omega,\, \cF)$ is said to be \embf{consistent}
with a $g$-function $g$ if $\P$ is shift-invariant and
\be{consist}
\int h(\omega)g(x\mid\omega)\, \P(d\omega) = \int_{\{\omega_0=x\}} h(\omega)\, \P(d\omega)
\ee
for all $x\in E$ and $\cF_{\leq -1}$-measurable function $h$. The family of these measures will be denoted
by $\mathcal{G}(g)$ and for each $\P\in\mathcal{G}(g)$, the process $(X_i)_{i\in\Z}$ on $(\Omega,\cF,\P)$
will be called a \embf{$g$-chain}.
\ed
\br{extrem}
In the consistency definition (\ref{consist}), $\P$ needs only to be defined on
$(\Omega_{-\infty}^{0},\, \cF_{\leq 0})$. Because of its shift-invariance, $\P$
can be extended in a unique way to $(\Omega,\, \cF)$. That's why, without loss
of generality, we can make no distinction between $\P$ on
$(\Omega_{-\infty}^{0},\, \cF_{\leq 0})$ and its natural extension on $(\Omega,\, \cF)$.
\er

\subsection{Regeneration}
\label{S1.4}

The following regeneration result is due to Comets, Fern\'andez and Ferrari \cite{comferfer02}
(see Theorem 4.1, Corollary 4.3 and Proposition 5.1). It will be the starting point of our
analysis.
\begin{theorem}[Comets \& al. (2002)]
\label{th0}
Let $g$ be a regular $g$-function such that
\be{hyp}
\prod_{k\geq 0} a_k >0
\quad\text{with}\quad
a_k=\inf_{\sigma_{-k}^{-1}\in\Omega_{-k}^{-1}} \, \sum_{\xi_0\in E}\,
\inf_{\omega_{-\infty}^{-k-1}\in \Omega_{-\infty}^{-k-1}}\,
g\bigl(\xi_0\mid\sigma_{-k}^{-1}\,\omega_{-\infty}^{-k-1}\,\bigr)
\ee
with the convention $\sigma_0^{-1}= \emptyset$.

Then
\begin{itemize}
\item[(i)] there exists a unique probability measure $\P$ consistent with $g$;
\item[(ii)] there exists a shift-invariant renewal process $(T_i)_{i\in\Z}$ with renewal distribution
\be{rwdist}
\P(T_{i+1}-T_i\geq M) = \rho_M,
\qquad M>0,\, i\neq 0
\ee
with $\rho_M$ the probability of return to the origin at epoch $M$ of the Markov chain on $\N\cup\{0\}$
starting at time zero at the origin with transition probabilities
\be{rwdist2}
\begin{cases}
p(k,k+1)=a_k,\\
p(k,0)=1-a_k,\\
p(k,j)=0 \text{ otherwise}
\end{cases}
\ee
and such that
\be{rwdist3}
T_0 \leq 0 < T_1.
\ee
\item[(iii)] the random blocks $\{(X_j \colon T_i\leq j <T_{i+1})\}_{i\in\Z}$, where $(X_i)_{i\in\Z}$
on $(\Omega,\cF,\P)$ is the associated $g$-chain, are independent and, except for $i=0$, identically
distributed.
\item[(iv)] $1\leq\E(T_2-T_1) = \sum_{i=1}^\infty \rho_i < \infty$.
\end{itemize}
\end{theorem}


\section{Main Results}
\label{S2}

We are now ready to state the main results of this paper.
\begin{theorem}[FCLT]
\label{th1}
Let $(X_i)_{i \in \Z}$ be a regular $g$-chain satisfying the regeneration assumption (\ref{hyp})
and $f\colon E\ra\R$ be a function such that
\be{th1eq2}
\E\big(f(X_0)\big)=0.
\ee
Then,
\be{th1eq3}
\frac{1}{\sqrt{n}} S_n
:= \frac{1}{\sqrt{n}}\sum_{i=0}^{n-1}f(X_i) \;\longrightarrow\; \cN\big(0,\sigma^2\big)
\quad \text{in distribution}
\ee
with
\be{var1}
0\,\leq\,
\sigma^2 = \frac{\E\left(\left(\sum_{i=T_1}^{T_2-1} f(X_i)\right)^2\right)}{\E(T_2-T_1)}
\,<\,\infty.
\ee
Furthermore, if $\E((T_2-T_1)^2)<\infty$, then
\be{var2}
\sigma^2=\E\Big(f^2(X_0)\Big)+2\sum_{i\geq1}\E\Big(f(X_0)\, f(X_i)\Big)<\infty.
\ee
\end{theorem}

\bc{cor1}
Let $(X_i)_{i \in \Z}$ be a regular $g$-chain such that there exist $a\in\R$, $b\in[1,\infty)$,
$C>0$ and $K\geq 2$ such that for all $k\geq K$
\be{cor1eq1}
a_k \geq 1-C\,\frac{(\log(k))^a}{k^b}.
\ee
Let $f\colon E\ra\R$ satisfying (\ref{th1eq2}).
\begin{itemize}
\item[(i)] If $a\in\R$ and  $b>1$, or $a<-1$ and $b=1$, then (\ref{th1eq3})
is satisfied with $\sigma^2$ defined by (\ref{var1}).
\item[(ii)] If $a\in\R$ and $b>2$, or $a<-1$ and $b=2$, then (\ref{th1eq3})
is satisfied with $\sigma^2$ defined by (\ref{var2}).
\end{itemize}
\ec
\br{rk1}
It is not easy to compare our results with previous ones based on mixing rates (see e.g.\ \cite{coh66},
\cite{ibrlin71} or \cite{griios90}), because there is no general relationship between mixing rates and
our continuity rates expressed in  terms of $(a_k)_{k\geq 0}$. However, in the case of regular $g$-chains
with exponential continuity decay, the mixing rates are also of exponential decay (see \cite{brefergal} or
\cite{fermai05}), and therefore CLT results as those describe in \cite{ibrlin71}, Chapter 18, are fulfilled.
\er

\section{Applications}
\label{S3}


\subsection{Binary autoregressive processes}
\label{S3.1}
The binary version of autoregressive processes is mainly used in statistics and econometrics.
It describes binary responses when covariates are historical values of the process (see
McCullagh and Nelder (1989), Section 4.3, for more details).

In what follows, we consider the example that was introduced previously in \cite{comferfer02}.
For the alphabet $E=\{-1,+1\}$, consider $\theta_0$ a real number and $(\theta_k:k\geq 1)$ an
absolutely summable real sequence. Let $q:\R \to (0,1)$ be a function strictly increasing and
continuously differentiable. Assume that $g(\,\cdot\mid\omega_{-\infty}^{-1})$ is the Bernoulli
law on $\{-1,+1\}$ with parameter $q(\theta_0+\sum_{k \geq 1} \theta_k \omega_{-k})$, that is,
\be{bap-1}
g\big(+1\mid\omega_{-\infty}^{-1}\big)
=q\bigg(\theta_0+\sum_{k \geq 1} \theta_k \omega_{-k}\bigg)
=1-g\big(-1\mid\omega_{-\infty}^{-1}\big).
\ee
Denote
\be{bap-3}
r_k=\sum_{m>k} |\theta_m|,
\qquad k\geq 0.
\ee

\bp{propAutoregress}
Let $(X_i)_{i \in \Z}$ be a regular $g$-chain with $g$ defined by (\ref{bap-1}) and
\be{bap-5}
f(x)=x-\E(X_0).
\ee
\begin{itemize}
\item[(i)] If $\sum_{k \geq 0} r_k^2 < \infty$, then $|\cG(g)|=1$.
\item[(ii)] If $\sum_{k \geq 0} r_k < \infty$, then the FCLT is satisfied with $\sigma^2$ defined
by (\ref{var1}) for $f$ defined by (\ref{bap-5}).
\item[(iii)] If there exists $K\geq 0$ such that $r_k \leq C(\log(k))^a/k^{b}$, $k\geq K$, with $C>0$, ($a\in\R$
and $b>2$) or ($a<-1$ and $b=2$), then the FCLT is satisfied with $\sigma^2$ defined by (\ref{var2}) for $f$
defined by (\ref{bap-5}).
\end{itemize}
\ep


\subsection{Power-law Ising chain}
\label{S3.2}

For the usual Ising (Gibbs) model, the central limit theorem is well known (see Newman \cite{new80}).
In the chain context, since the relationship between one dimensional Gibbs measures and chains, discussed in
\cite{fermai04}, does not allow to easily interpret Gibbs results in a chain setting, the problem becomes relevant.
In what follows, we give central limit theorem for power-law Ising chains.

For the alphabet $E=\{-1,1\}$, consider the power-law Ising chain defined by
\be{plIC-def}
g\big(\omega_0\mid\omega_{-\infty}^{-1}\big)
=\frac{\displaystyle\exp\Bigg[-\sum_{k=-\infty}^{-1}\phi_k\big(\omega_{-\infty}^{0}\big)\Bigg]}
{\displaystyle\sum_{\sigma_0\in E}\exp\Bigg[-\sum_{k=-\infty}^{-1}\phi_k\big(\omega_{-\infty}^{-1}\sigma_0\big)\Bigg]}
\ee
with
\be{plIC-def*}
\phi_k\big(\omega_{-\infty}^{0}\big)=-\beta\frac{1}{|k|^p}\, \omega_0\, \omega_k,
\qquad k\leq -1,\beta>0,\, p>1.
\ee

\bp{propIsing}
Let $(X_i)_{i \in \Z}$ be a regular $g$-chain with $g$ defined by (\ref{plIC-def}--\ref{plIC-def*}) and
\be{plIC-5}
f(x)=x.
\ee
\begin{itemize}
\item[(i)] If $p>3/2$, then $|\cG(g)|=1$.
\item[(ii)] If $p>2$, then the FCLT is satisfied with $\sigma^2$ defined by (\ref{var1})
for $f$ defined by (\ref{plIC-5}).
\item[(iii)] If $p>3$, then the FCLT is satisfied with $\sigma^2$ defined by (\ref{var2})
for $f$ defined by (\ref{plIC-5}).
\end{itemize}
\ep


\section{Proofs}
\label{S4}


\subsection{Proof of Theorem \ref{th1}}
\label{S4.1}

Take $(T_i)_{i \in \Z}$ as given by Theorem \ref{th0} (ii) and  write
\be{p1-9}
S_n=\sum_{i=0}^{n-1} f(X_i)=\sum_{i=0}^{T_1-1} f(X_i)
+ \sum_{i=T_1}^{T_{i(n)}-1} f(X_i) + \sum_{i=T_{i(n)}}^{n-1} f(X_i)
\ee
with
\be{p1-11}
i(n)=\left\{\ba{ll}\max\{k \geq 1\colon T_k < n\} & \text{ if } T_1 < n,\\0 & \text{ otherwise.}\ea\right.
\ee
\bl{zeroprob-lem}
Under (\ref{hyp}), both $n^{-1/2}\sum_{i=0}^{T_1-1} f(X_i)$ and $n^{-1/2}\sum_{i=T_{i(n)}}^{n-1} f(X_i)$
tend to zero in probability.
\el
\bpr
For any $K>0$, we obviously have
\be{p1-13}
\begin{aligned}
\P\Bigg(\bigg|\sum_{i=0}^{T_1-1} f(X_i)\bigg|>K\sqrt{n}\Bigg)
&\leq \P\left(T_1 \sup_{x \in E} \left|f(x)\right|>K\sqrt{n}\right).
\end{aligned}
\ee
Therefore,
\be{p1-14}
\sup_{n \geq 1}\, \P\Bigg(\bigg|\sum_{i=0}^{T_1-1} f(X_i)\bigg|>K\sqrt{n}\Bigg) \leq \P(T_1 > M_K)
\ee
with
\be{MKdef}
M_K=\frac{K}{M}
\quad\text{and}\quad
M=\sup_{x\in E} \left|f(x)\right|.
\ee
On the other hand, by the shift-invariance of $(T_i)_{i\in\Z}$,
\be{p1-15}
\begin{aligned}
\P\Bigg(\bigg|\sum_{i=T_{i(n)}}^{n-1} f(X_i)\bigg|>K\sqrt{n}\Bigg)
&\leq \P\Bigg(\sum_{i=T_{i(n)}}^{T_{i(n)+1}-1} \left|f(X_i)\right|>K\sqrt{n}\Bigg)\\
&\leq \P\left((T_{i(n)+1}-T_{i(n)}) \sup_{x\in E} \left|f(x)\right|>K\sqrt{n}\right)\\
&= \P\left((T_1-T_0) \sup_{x\in E} \left|f(x)\right|>K\sqrt{n}\right),
\end{aligned}
\ee
where we used that $i(0)=0$. Therefore,
\be{p1-16}
\sup_{n \geq 1}\, \P\Bigg(\bigg|\sum_{i=T_{i(n)}}^{n-1} f(X_i)\bigg|>K\sqrt{n}\Bigg)
\leq \P\big(T_1-T_0 > M_K\big),
\ee
where $M_K$ is defined by (\ref{MKdef}). Noticing that, under condition (\ref{hyp}), both $\P(T_1 > M_K)$
and $\P(T_1-T_0 > M_K)$ tend to zero as $K$ goes to infinity, we can conclude that both
$n^{-1/2}\sum_{i=0}^{T_1-1} f(X_i)$ and $n^{-1/2}\sum_{i=T_{i(n)}}^{n-1} f(X_i)$
tend to zero in probability.
\epr

\bp{prop1}
Let $(X_i)_{i \in \Z}$ be a regular $g$-chain satisfying (\ref{hyp}) and $f\colon E\ra\R$.
Then the following statements are equivalent:
\begin{itemize}
\item[(i)] $S_n/\sqrt{n} \;\longrightarrow\; \cN(0,\sigma^2)$ in distribution for
\be{prop1eq1}
\sigma^2 = \frac{\E\left(\left(\sum_{i=T_1}^{T_2-1} f(X_i)\right)^2\right)}{\E(T_2-T_1)},
\ee
where $T_1$ and $T_2$ are defined in(\ref{rwdist}--\ref{rwdist3});
\item[(ii)] $(S_n/\sqrt{n})_{n\geq 0}$ is bounded in probability;
\item[(iii)] $\E\Big(\sum_{i=T_1}^{T_2-1}f(X_i)\Big)=0$
and $\E\bigg(\Big(\sum_{i=T_1}^{T_2-1}f(X_i)\Big)^2\bigg)<\infty$.
\end{itemize}
\ep

\bpr
The direction (i) $\Rightarrow$ (ii) is trivial, so we only have to show (ii) $\Rightarrow$ (iii) and
(iii) $\Rightarrow$ (i).

To prove (ii)$\Rightarrow$ (iii), we first remark that equation (\ref{p1-9}), Lemma \ref{zeroprob-lem}
and assertion (ii) imply that $n^{-1/2}\sum_{i=T_1}^{T_{i(n)}-1} f(X_i)$ is bounded in probability.
Then, by the converse of the central limit theorem for real i.i.d.\ sequences (see e.g.\ \cite{ledtal91},
Section 10.1), we must have (iii).

To prove (iii) $\Rightarrow$ (i), we first see that
\be{p1-19}
\frac{i(n)}{n} \;\longrightarrow\; \frac{1}{\E(T_2-T_1)} \quad \text{a.s.},
\ee
which follows from Theorem 5.5.2 of \cite{chu74} and the fact that $(T_{i+1}-T_i)_{i>0}$ is
an i.i.d.\ process. Let us denote
\be{p1-21}
\xi_k=\sum_{i=T_k}^{T_{k+1}-1} f(X_i)
\quad \text{and} \quad
e(n)=\left\lfloor\frac{n}{\E(T_2-T_1)}\right\rfloor.
\ee
Thanks to the Lemma \ref{zeroprob-lem}, to prove (i), it is enough to show that
\be{p1-45}
\frac{1}{\sqrt{n}}\sum_{k=1}^{i(n)-1} \xi_k
\;\longrightarrow\;
\mathcal{N}\left(0,\frac{\E(\xi_1^2)}{\E(T_2-T_1)}\right) \quad \text{in law.}
\ee
This follows from the standard central limit theorem result
\be{p1-31}
\frac{1}{\sqrt{n}}\sum_{k=1}^{e(n)} \xi_k
\;\longrightarrow\;
\mathcal{N}\left(0,\frac{\E(\xi_1^2)}{\E(T_2-T_1)}\right) \quad \text{in law,}
\ee
and
\be{p1-29}
\frac{1}{\sqrt{n}}\Bigg(\sum_{k=1}^{i(n)-1} \xi_k - \sum_{k=1}^{e(n)} \xi_k\Bigg)
\;\longrightarrow\; 0 \quad \text{in probability.}
\ee
To prove the latter, for any $\epsilon>0$, let
\be{p1-23}
K'=\left\lfloor\frac{\epsilon^3}{2\big(1+\E(\xi_k^2)\big)}\right\rfloor.
\ee
First, remark that
\be{p1-24}
\begin{aligned}
&\P\Bigg(\bigg|\sum_{k=1}^{i(n)-1} \xi_k - \sum_{k=1}^{e(n)} \xi_k\bigg| \geq \sqrt{n}\epsilon\Bigg)\\
&\qquad\leq 2\, \P\Bigg(\max_{j\leq nK'}\bigg|\sum_{k=1}^j \xi_k\bigg| \geq \sqrt{n}\epsilon\Bigg)
+\P\bigg(\big|i(n)-1-e(n)\big| \geq nK'\bigg).
\end{aligned}
\ee
Then, use the Kolmogorov's maximal inequality, (iii) and (\ref{p1-23}), to get
\be{p1-25}
\P\Bigg(\max_{j\leq nK'}\Bigg|\sum_{k=1}^j \xi_k\Bigg| \geq \sqrt{n}\epsilon\Bigg)
\leq \frac{1}{n\epsilon^2}\; \Var\Bigg(\sum_{k=1}^{nK'} \xi_k\Bigg)
= \frac{K'\;\E(\xi_k^2)}{\epsilon^2}
\leq \frac{\epsilon}{2}.
\ee
Since
\be{p1-27}
\lim_{n \to \infty} \P\bigg(\big|i(n)-1-e(n)\big| \geq nK'\bigg) = 0,
\ee
it follows from (\ref{p1-24}--\ref{p1-25}) that
\be{p1-28}
\lim_{n\to\infty}
\P\Bigg(\bigg|\sum_{k=1}^{i(n)-1} \xi_k - \sum_{k=1}^{e(n)}\xi_k\bigg|\geq\sqrt{n}\epsilon\Bigg)
<\epsilon.
\ee
\end{proof}

\bprthm{th1}
First, using the shift-invariance of $\P$, we see that for all $n>1$
\be{p2-3}
\begin{aligned}
\frac{\E(S_n^2)}{n} &= \frac{1}{n}\Bigg(\E\Bigg(\sum_{i=0}^{n-1} f^2(X_i)\Bigg)
+2\,\E\Bigg(\sum_{i=0}^{n-2} \sum_{j=i+1}^{n-1} f(X_i)f(X_j)\Bigg)\Bigg)\\
&= \E\Big(f^2(X_0)\Big)+\frac{2}{n}\sum_{i=1}^{n-1} (n-i)\,\E\Big(f(X_0)f(X_i)\Big).
\end{aligned}
\ee
Then, under assumptions (\ref{hyp}) and (\ref{th1eq2}), Theorem 1 of \cite{brefergal} insures that
\be{p2-4}
\sum_{i\geq 1}\E\Big(f(X_0)f(X_i)\Big)
\leq C\sum_{i\geq 1}\rho_i<\infty,
\ee
for some $C>0$. Therefore, it follows from Kronecker's lemma that
\be{p2-5}
\lim_{n \to \infty} \frac{\E(S_n^2)}{n} = \E\Big(f^2(X_0)\Big)
+2\sum_{i\geq 1}\E\Big(f(X_0)f(X_i)\Big)<\infty.
\ee
In particular, Chebyshev's inequality implies that
\be{p2-7}
\left(\frac{S_n}{\sqrt{n}}\right)_{n \geq 1}
\ee
is bounded in probability. Then Proposition \ref{prop1} concludes the proof of (\ref{th1eq3}-\ref{var1}).

To show (\ref{var2}), we use first Proposition \ref{prop1} (iii) to get
\be{p2-9}
\E\Bigg(\bigg[\sum_{k=1}^{i(n)+1} \xi_k \bigg]^2\Bigg)
= \E\Bigg(\sum_{k=1}^{i(n)+1} \xi_k^2\Bigg)
=\E\big(\xi_k^2\big)\, \E\big(i(n)+1\big),
\ee
where the rightmost equality follows from the Wald's Equation and the fact that $i(n)+1$ is a stopping time
w.r.t.\ $(\xi_k)_{k\geq 1}$. Since
\be{p2-10}
\lim_{n\to\infty}\frac{\E\big(i(n)\big)}{n}=\frac1{\E(T_2-T_1)}
\ee
(see Theorem 5.5.2 in \cite{chu74}), (\ref{var1}) and (\ref{p2-9}) give
\be{p2-10*}
\lim_{n\to\infty}\frac1n\, \E\Bigg(\bigg[\sum_{k=1}^{i(n)+1} \xi_k \bigg]^2\Bigg)
=\sigma^2.
\ee
But, we have
\be{p2-11}
\lim_{n \to \infty} \frac{1}{n} \, \E\left(\max_{k \leq n} \xi_k^2\right)=0
\ee
(see \cite{chu67} p.90 for a proof) and therefore
\be{p2-12}
\lim_{n \to \infty} \frac{1}{n} \, \E\big(\xi_{k}^2\big)=0,
\qquad k\in\{i(n), i(n)+1\}.
\ee
Hence, by (\ref{p2-10*}) and (\ref{p2-12}), we finally have
\be{p2-13}
\lim_{n \to \infty} \frac{1}{n} \, \E\Bigg(\bigg[\sum_{k=1}^{i(n)-1} \xi_k\bigg]^2\Bigg)=\sigma^2.
\ee
Now, if we show that
\be{p2-15}
\lim_{n \to \infty} \frac{1}{n} \, \E\Bigg(\bigg[\sum_{k=1}^{T_1-1}f(X_k)\bigg]^2\Bigg)=0
\ee
and
\be{p2-17}
\lim_{n \to \infty} \frac{1}{n} \, \E\Bigg(\bigg[\sum_{k=T_{i(n)}}^{n-1}f(X_k)\bigg]^2\Bigg)=0,
\ee
then, using Theorem \ref{th0} (iii), (\ref{p1-9}) and (\ref{p2-13}), we can conclude the proof
of (\ref{var2}).

The proofs of (\ref{p2-15}--\ref{p2-17}) are similar, we will prove (\ref{p2-17}) only. To that aim, we
first remark that
\be{p2-19}
k\, \P(T_1-T_0\geq k) \rightarrow 0 \quad \text{as} \quad k \to \infty,
\ee
because $\P(T_1-T_0\geq k)$ is decreasing to 0 and
\be{p2-20}
\sum_{k \geq 0} \P(T_1-T_0\geq k)
= \E(T_1-T_0)
=\frac{\E([T_2-T_1]^2)}{\E(T_2-T_1)}< \infty,
\ee
where the rightmost equality uses that $\P(T_1-T_0=k)=k\P(T_2-T_1=k)/\E(T_2-T_1)$
(see e.g.\ Lawler \cite{law06}, Chapter 6, Section 6.2, equality (6.10)). Therefore, recalling (\ref{MKdef})
and using the shift-invariance of $\P$, we have
\be{p2-21}
\begin{aligned}
\frac{1}{n} \, \E\Bigg(\Bigg[\sum_{k=T_{i(n)}}^{n-1}f(X_k)\Bigg]^2\Bigg) &\leq \frac{M^2}{n} \,
\E\Big([n-T_{i(n)}]^2\Big)\\
&= \frac{M^2}{n} \, \sum_{k=1}^n k^2 \, \P\big(n-T_{i(n)} = k\big)\\
&= \frac{M^2}{n} \, \sum_{k=1}^n (2k-1) \P\big(n-T_{i(n)} \geq k\big)\\
&\leq \frac{M^2}{n} \, \sum_{k=1}^n (2k-1)\, \P(T_1-T_0\geq k),
\end{aligned}
\ee
which in view of (\ref{p2-19}), goes to zero as $n$ tends to infinity.
\eprthm


\subsection{Proof of Corollary \ref{cor1}}
\label{S4.2}

\bpr
To prove (i), it suffices to see that when $(a_k)_{k\geq 0}$ satisfies (\ref{cor1eq1}),
then ($a\in\R$ and $b>1$) or ($a<-1$ and $b=1$) if and only if $\sum_{k\geq 0}(1-a_k)<\infty$
which is equivalent to $\prod_{k\geq 0}a_k>0$. Therefore, applying Theorem \ref{th1} first part,
we get the result.

To prove (ii), we first denote
\be{aktilde}
\tilde{a}_k = 1-C\,\frac{(\log(k))^a}{k^b}
\ee
and the associated $\tilde\rho_k$, defined by the analogous of (\ref{rwdist2}). Since $b>1$,
Proposition 5.5 (iv) in Fern\'andez, Ferrari and Galves \cite{ferfergal01} gives that there
exists some constant $C_1>0$ so that
\be{}
\tilde\rho_k\leq C_1(1-\tilde a_k),
\qquad k\geq 0.
\ee
Therefore, using that $a_k\geq \tilde a_k$ implies $\rho_k\leq \tilde \rho_k$, we have
\be{}
\begin{aligned}
\E\big((T_2-T_1)^2\big)
=\sum_{k\geq 0} k\rho_k
\leq C_1\sum_{k\geq 0} k(1-\tilde{a}_k),
\end{aligned}
\ee
which is finite if and only if ($a\in\R$ and $b>2$) or ($a<-1$ and $b=2$).
Then, applying Theorem \ref{th1} second part, we get the result.
\epr


\subsection{Proof of Proposition \ref{propAutoregress}}
\label{S4.3}

\bpr
Define the variation by
\be{pa-7}
\var_k
=\sup\Big\{\big|g\big(\omega_0\mid\omega_{-\infty}^{-1}\big)
-g\big(\sigma_0\mid\sigma_{-\infty}^{-1}\big)\big|
\colon \omega_{-\infty}^0,\sigma_{-\infty}^0 \in \Omega_{-\infty}^0,\,
\omega_{-k}^0=\sigma_{-k}^0\Big\},
\quad k\geq 0.
\ee
Then, because $|E|=2$, we have for any $k\geq 0$ (with the convention $\sigma_{0}^{-1}=\emptyset$)
\be{akvark}
\begin{aligned}
a_k
&=\inf\Big\{
g\big(1\mid \omega_{-\infty}^{-k-1}\,\sigma_{-k}^{-1}\big)
+g\big(-1\mid \xi_{-\infty}^{-k-1}\,\sigma_{-k}^{-1}\big)\colon
\sigma_{-k}^{-1}\in\Omega_{-k}^{-1},\, \omega_{-\infty}^{-k-1},\xi_{-\infty}^{-k-1}\in\Omega_{-\infty}^{-k-1}\Big\}\\
&=1-\sup\Big\{
-g\big(1\mid \omega_{-\infty}^{-k-1}\,\sigma_{-k}^{-1}\big)
+g\big(1\mid \xi_{-\infty}^{-k-1}\,\sigma_{-k}^{-1}\big)\colon
\sigma_{-k}^{-1}\in\Omega_{-k}^{-1},\, \omega_{-\infty}^{-k-1},\xi_{-\infty}^{-k-1}\in\Omega_{-\infty}^{-k-1}\Big\}\\
&=1-\var_k.
\end{aligned}
\ee
Therefore,
\be{pa-9}
a_k=1-\sup\Bigg\{ q\Bigg(\theta_0+\sum_{j=1}{k} \theta_j\, \sigma_{-j}+r_k\Bigg)
-q\Bigg(\theta_0+\sum_{j=1}^{k} \theta_j\, \sigma_{-j}-r_k\Bigg)\colon
\sigma_{-k}^{-1}\in\Omega_{-k}^{-1}\Bigg\}.
\ee
Because $q$ is continuously differentiable on a compact set, there exists $C>0$ such that
\be{pa-11}
\var_k \leq C r_k,
\qquad k \geq 0.
\ee

To prove (i), it suffices to remark that
\be{Obejoh}
\sum_{k\geq 0}\var_k^2<\infty,
\ee
which is a tight uniqueness criteria in terms on the variation (see Johansson and \"Oberg \cite{johobe03}
and Berger, Hoffman and Sidoravicius \cite{berhofsid05}).

To show (ii), we simply note that
\be{pa-13}
\E\big(f(X_0)\big)=\E\big(X_0-\E(X_0)\big)=0
\ee
and that
\be{pa-15}
\sum_{k \geq 0} r_k < \infty \quad\Longrightarrow\quad \prod_{k \geq 0} a_k > 0 .
\ee

Finally, to prove part (iii), it suffices first to combine (\ref{akvark}) and (\ref{pa-11}) to get
\be{}
a_k\geq 1-Cr_k,
\qquad k\geq 0,
\ee
and then to apply Corollary \ref{cor1} (ii).

\epr


\subsection{Proof of Proposition \ref{propIsing}}
\label{S4.4}

\bpr
We need the following well-known bound, whose proof we present for completeness is given
in Appendix \ref{A} (we follow the approach of~\cite[Lemma V.1.4]{sim93}).
\bl{simlem}
Let $g$ be a $g$-function satisfying (\ref{plIC-def}) with $\sum_{k=-\infty}^{-1}|\phi_k|<\infty$
and $h$ be a $\cF_0$-measurable function. Then, for any $\omega_{-\infty}^{-1}$,
$\sigma_{-\infty}^{-1}\in\Omega_{-\infty}^{-1}$,
\be{sim-ineq}
\begin{aligned}
&\bigg|\sum_{\omega_0\in E}h(\omega_0)\Big(g\big(\omega_0\mid\omega_{-\infty}^{-1}\big)
-g\big(\omega_0\mid\sigma_{-\infty}^{-1}\big)\Big)\bigg|\\
&\qquad\leq \sup_{x\in E}|h(x)|\, \sup_{\omega_0\in E}\bigg|\sum_{k=-\infty}^{-1}
\Big[\phi_k\big(\omega_{-\infty}^{0}\big)-\phi_k\big(\sigma_{-\infty}^{-1}\omega_0\big)\Big]\bigg|.
\end{aligned}
\ee
\el
Applying the previous lemma for $h\equiv1$, we obtain that, for any
$\omega_{\le i-1}, \sigma_{\le i-1}\in \Omega_{\le i-1}$,
\be{Ising-3}
\Big|g\big(\omega_0\mid\omega_{-\infty}^{-1}\big)-g\big(\omega_0\mid\sigma_{-\infty}^{-1}\big)\Big|
\leq |\beta| \sup_{\omega_0\in E}\bigg|\sum_{k=-\infty}^{-1}\frac{1}{|k|^p}\big(\omega_0\omega_k-
\omega_0\sigma_k\big)\bigg|,
\ee
from which
\be{Ising-5}
\var_k
\leq 2|\beta|\sum_{j=-\infty}^{-k-1}\frac{1}{|j|^p}
\leq 2|\beta|\frac{1}{k^{p-1}},
\qquad k\geq 1,
\ee
is an immediate consequence.

To prove (i), it suffices to see that (\ref{Ising-5}) with $p>3/2$ implies the validity of
(\ref{Obejoh}).

To prove (ii) and (iii), we first remark that similarly to the gibbsian setting, it can be easily checked
that $\E(X_0)=0$ and therefore under (\ref{plIC-5}), (\ref{th1eq2}) is fulfilled. Then, we combine
(\ref{akvark}) and (\ref{Ising-5}) to get
\be{Ising-7}
a_k \geq 1-2|\beta|\frac{1}{k^{p-1}},
\quad k\geq 1.
\ee
Thus, the results are direct consequences of Corollary \ref{cor1} (i) and (ii).
\epr


\appendix

\section{Appendix}
\label{A}

In this appendix we give the proof of Lemma \ref{simlem}.

\bpr
For all $\omega_{-\infty}^{-1},\sigma_{-\infty}^{-1}\in\Omega_{-\infty}^{-1}$ and $0<\theta<1$, define
$\Gamma_{\omega,\sigma}^{\theta}\colon E\ra(0,1)$ by
\be{simlem-11}
\Gamma_{\omega,\sigma}^{\theta}(\xi_0)
=\frac{\displaystyle\exp\Big[\theta H_{\omega}(\xi_0)+(1-\theta)H_{\sigma}(\xi_0)\Big]}
{\displaystyle\sum_{\eta_0\in E}\exp\Big[\theta H_{\omega}(\eta_0)+(1-\theta)H_{\sigma}(\eta_0)\Big]}
\quad\text{with}\quad
H_\omega(\xi_0)=\sum_{k=-\infty}^{-1}\phi_k\big(\omega_{-\infty}^{-1}\xi_0\big).
\ee
Then, to prove (\ref{sim-ineq}), it suffices to see that
\be{simlem-13}
\begin{aligned}
&\bigg|\sum_{\xi_0\in E}h(\xi_0)\Big(g\big(\xi_0\mid\omega_{-\infty}^{-1}\big)
-g\big(\xi_0\mid\sigma_{-\infty}^{-1}\big)\Big)\bigg|\\
&\quad\leq \int_{0}^{1}\bigg|\frac{d}{d\theta}\bigg[\sum_{\xi_0\in E} h(\xi_0)\,
\Gamma_{\omega,\sigma}^{\theta}(\xi_0)\bigg]\bigg|\, d\theta\\
&\quad=\int_{0}^{1}\bigg|\sum_{\xi_0\in E} h(\xi_0)\Big(H_{\omega}(\xi_0)
-H_{\sigma}(\xi_0)\Big) \Gamma_{\omega,\sigma}^{\theta}(\xi_0)\\
&\qquad-\sum_{\xi_0\in E} h(\xi_0)\, \Gamma_{\omega,\sigma}^{\theta}(\xi_0)
\sum_{\eta_0\in E} \Big(H_{\omega}(\eta_0)-H_{\sigma}(\eta_0)\Big)
\Gamma_{\omega,\sigma}^{\theta}(\eta_0)\bigg|\, d\theta\\
&\quad\leq \|h\|_\infty\int_{0}^{1}\sum_{\xi_0\in E} \bigg|\Big(H_{\omega}(\xi_0)
-H_{\sigma}(\xi_0)\Big) -\sum_{\eta_0\in E} \Big(H_{\omega}(\eta_0)
-H_{\sigma}(\eta_0)\Big)\Gamma_{\omega,\sigma}^{\theta}(\eta_0) \bigg|\,
\Gamma_{\omega,\sigma}^{\theta}(\xi_0)\, d\theta\\
&\quad\leq \|h\|_\infty\sup_{\xi_0 E}\big|H_{\omega}(\xi_0)
-H_{\sigma}(\xi_0)\big|.
\end{aligned}
\ee
\epr


\end{document}